\documentstyle[twoside,amstex]{article}
\input amssym.def
\input amssym.tex
\def\bc{\begin{center}}
\def\ec{\end{center}}
\def\no{\noindent}
\def\hang{\hangindent\parindent}
\def\textindent#1{\indent\llap{[#1]\enspace}\ignorespaces}
\def\re{\par\hang\textindent}
\begin{document}
\thispagestyle{empty} \vspace*{3 true cm} \pagestyle{myheadings}
\markboth {\hfill {\sl M.S. Abdolyousefi and H. Chen}\hfill}
{\hfill{\sl Elementary Matrix Reduction Over Certain Rings}\hfill} \vspace*{-1.5 true cm} \bc{\large\bf Elementary Matrix Reduction Over Certain Rings}\ec

\vskip6mm
\bc{{\bf Marjan
Sheibani Abdolyousefi}\\[2mm]
Faculty of Mathematics, Statistics and Computer Science\\
Semnan University, Semnan, Iran\\
m.sheibani1@@gmail.com}\ec

\vskip2mm
\bc{{\bf Huanyin Chen$^*$}\\[2mm]
Department of Mathematics, Hangzhou Normal University\\
Hangzhou 310036, China\\
huanyinchen@@aliyun.com}\ec

\begin{figure}[b]
\vspace{-3mm}
\rule[-2.5truemm]{5cm}{0.1truemm}\\[2mm]
{\footnotesize 2010 Mathematics Subject Classification. 13F99, 13E15,
16E50.
$^*$ Corresponding author.}
\end{figure}

\vskip10mm
\begin{minipage}{120mm}
\no {\bf Abstract:} We explore elementary matrix reduction over certain rings characterized by their localizations.
Let $R$ be a locally stable ring. We prove that $R$ is an elementary divisor ring if and only if it is a B$\acute{e}$zout ring.
For any matrix $A$ of size $n$ we further prove that there exist unimodular matrices $P\in GE_n(R)$ and $Q\in GL_n(R)$, that $PAQ=diag(\varepsilon_1,\cdots ,\varepsilon_n)$ where $\varepsilon_i$ is a divisor of $\varepsilon_{i+1}$.
Elementary reduction of matrices over some related localizations are also studied. Many known results on domains are thereby extended to general commutative rings which may contain many zero divisors.
\end{minipage}

\vskip15mm \bc{\bf 1. Introduction}\ec

\vskip4mm \no Throughout this paper, all rings are commutative
with an identity. A matrix $A$ (not necessarily square) over a ring
$R$ admits diagonal reduction if there exist invertible matrices
$P$ and $Q$ such that $PAQ$ is a diagonal matrix $(\varepsilon_{ij})$, for
which $\varepsilon_{ii}$ is a divisor of $\varepsilon_{(i+1)(i+1)}$ for each $i$. A
ring $R$ is called an elementary divisor ring provided that every
matrix over $R$ admits a diagonal reduction. A ring $R$ is a B$\acute{e}$zout ring if every finitely generated ideal is principal. It is well known that
every elementary divisor ring is a B$\acute{e}$zout ring. But the converse is not true. An attractive problem is to investigate various conditions under which a B$\acute{e}$zout ring is an elementary divisor ring. After Kaplansky's work on
elementary divisor rings without nonzero zero divisors, Gillman and Henriksen proved that

\vskip2mm \hspace{-1.8em} {\bf Theorem 1.1 [4, Theorem 1.1].}\ \
{\it A ring $R$ is an elementary divisor ring if and only if}
\begin{enumerate}
\item [(1)] {\it $R$ is a Hermite ring;}
\vspace{-.5mm}
\item [(2)] {\it For all $a_1,a_2,a_3\in R$, $a_1R+a_{2}R+a_{3}R=R\Longrightarrow $ there exist $p,q\in R$ such that
$(pa_1+qa_2)R+qa_{3})R=R$.}\end{enumerate}

Here, a ring $R$ is a Hermite ring if and only if for any $a,b\in R$ there exists some $Q\in GL_2(R)$ such that $(a,b)Q=(c,0)$ for some $c\in R$. As is well known, a domain $R$ (i.e., a ring without zero divisors) is a Hermite ring if and only if $R$ is a B$\acute{e}$zout ring (cf. [7, Theorem 1.2.6]). But they do not coincide with each other for rings with non-zero zero divisors. Let $S$ be a multiplicatively closed set of $R$, i.e., $1$ is in $S$ and for $s$ and $t$ in $S$ we also have $st$ in $S$. We use $R_S$ to denote the localization which is a ring consisting of "fractions" with denominators coming from $S$.
The motivation of this paper is to explore elementary matrix reduction over certain rings characterized by their localizations.

Let $a\in R$. We denote $S_a$ as a set $\{ b\in R~\mid~aR+bR=R\}$.
Then $S_a$ is a multiplicatively closed set. Let $R_a=S_a^{-1}R$. We say that $R_a$ is a localized ring of $R$ at $a\in R$.
A ring $R$ has stable range 1 provided that $cR+dR=R$ with $c,d\in R\Longrightarrow ~\exists y\in R$ such that $c+dy\in U(R)$.
This condition plays an important rule in algebraic $K$-theory. The class of rings having stable range 1 is very large. It includes exchange rings, regular rings, $\pi$-regular rings, semilocal rings, local rings, etc. For general theory of such condition, we refer the read to [1]. We say that a ring $R$ is locally stable if $cR+dR=R$ with $c,d\in R\Longrightarrow ~\exists y\in R$ such that $R_{c+dy}$ has stable range 1. We shall see that many known rings belong such type of rings.

Let $R$ be a locally stable ring. We prove that $R$ is an elementary divisor ring if and only if it is a B$\acute{e}$zout ring.
For any matrix $A$ of size $n$ we further prove that there exist unimodular matrices $P\in GE_n(R)$ and $Q\in GL_n(R)$, that $PAQ=diag(\varepsilon_1,\cdots ,\varepsilon_n)$ where $\varepsilon_i$ is a divisor of $\varepsilon_{i+1}$. Elementary reduction of matrices over a B$\acute{e}$zout ring $R$ such that $R_a$ or $R_{1-a}$ has stable range 1 for all $a\in R$ is studied as well. An element $c\in R$ is adequate provided that for any $a\in R$ there exist some $r,s\in R$ such that $(1)$ $c=rs$; $(2)$ $rR+aR=R$; $(3)$ $s'R+cR\neq R$ for each non-invertible divisor $s'$ of $s$.
Following Domsha and Vasiunyk,, a ring $R$ has adequate range 1 provided that $aR+bR=R$ with $a,b\in R\Longrightarrow ~\exists y\in R$ such that $a+by\in R$ is adequate [2]. For instance, rings having stable range 1, adequate rings and VNL rings have adequate range 1. Let $T=\{ a~|~R_a~$ has stable range one $\}$ and $S=\{ a\in R~|~a\in R$ is adequate $\}$. Then $T$ and $S$ are multiplicative closed sets of $R$. Finally, we prove that if $R_T$ or $R_S$ has adequate range 1, then $R$ is an elementary divisor ring if and only if $R$ is a Hermite ring. Many known results are thereby
generalized to much wider class of rings, e.g. [2, Theorem 14], [4, Theorem 3.7], [7, Theorem 3.2.2, Theorem 4.3.5 and Theorem 4.9.1], [8, Proposition 10] and [10, Proposition 5].

We shall use $J(R)$ and $U(R)$ to denote the Jacobson radical of $R$ and the set of all units in $R$, respectively. $M_n(R)$ denotes
the ring of all $n\times n$ matrices over $R$, and $GL_n(R)$ stands for the $n$-dimensional general linear group of $R$. Denote the subgroup of $GL_n(R)$ generated by elementary matrices by $GE_n(R)$.

\vskip10mm\bc{\bf 2. Locally Stable Rings}\ec

\vskip4mm The purpose of this section is to investigate elementary properties of locally stable rings, and then explore elementary reduction of matrices over such rings. We begin with

\vskip4mm \hspace{-1.8em} {\bf Lemma 2.1.}\ \ {\it Let $R$ be a ring, and let $a\in R$. Then $a\in J(R_a)$.}
\vskip2mm\hspace{-1.8em} {\it Proof.}\ \ For any $\frac{r}{s}\in R_a$, we check that $\frac{1}{1}-\frac{a}{1}\cdot \frac{r}{s}=\frac{s-ar}{s}$. Clearly, $s\in S_a$, i.e., $sR+aR=R$. Write $sp+aq=1$, and so $(s-ar)p+a(q+rp)=1$. It follows that $(s-ar)R+aR=R$, and so $\frac{s-ar}{s}\in U(R_a)$. Hence, $\frac{a}{1}\in J(R_a)$. This completes the proof.\hfill$\Box$

\vskip4mm \hspace{-1.8em} {\bf Lemma 2.2.}\ \ {\it Let $R$ be a ring and let $0\neq a\in R$. Then the following are equivalent:}
\vspace{-.5mm}
\begin{enumerate}
\item [(1)]{\it $R_a$ has stable range 1.}
\vspace{-.5mm}
\item [(2)]{\it $R/aR$ has stable range 1.}
\end{enumerate}
\vspace{-.5mm} {\it Proof.}\ \  $(1)\Rightarrow (2)$ Given $\overline{b}(R/aR)+\overline{c}(R/aR)=R/aR$, then we have some $x,y\in R$ such that $(bx+cy)R+aR=R$. This implies that $bx+cy\in S_a$. Hence, $\frac{b}{1}R_a+\frac{c}{1}R_a=R_a$. By hypothesis, we have a $\frac{z_1}{z_2}\in R_a$ such that
$\frac{b}{1}+\frac{c}{1}\cdot \frac{z_1}{z_2}\in U(R_a)$. This implies that $\frac{bz_2+cz_1}{z_2}\in U(R_a)$, and so $bz_2+cz_1\in S_a$. Hence, $(bz_2+cz_1)R+aR=R$.
As $z_2\in S_a$, we see that $z_2R+aR=R$, and so $\overline{z_2}\in U(R/aR)$. We infer that $\overline{b}+\overline{cz_1}\cdot \overline{z_2}^{-1}\in U(R/aR)$. Therefore $R/aR$ has stable range 1.

$(2)\Rightarrow (1)$ Let $\frac{b}{s}R_a+\frac{c}{s}R_a=R_a$. Then $\frac{bb'}{ss'}+\frac{cc'}{st'}=\frac{1}{1}$. Hence, we can find some $u,v\in R$ such that $bu+cv=ss't'\in S_a$. Thus, we get $aR+(bu+cv)R=R$, and so $\overline{b}(R/aR)+\overline{c}(R/aR)=R/aR$. By hypothesis, there exists a $y\in R$ such that $\overline{b+cy}\in U(R/aR)$. This shows that $(b+cy)R+aR=R$; hence that $\frac{b}{s}+\frac{c}{s}\cdot \frac{y}{1}=\frac{b+cy}{s}\in U(R_a)$. Therefore $R_a$ has stable range 1, as asserted.\hfill$\Box$

\vskip4mm \hspace{-1.8em} {\bf Theorem 2.3.}\ \ {\it Let $R$ be a ring. Then the following are equivalent:}
\vspace{-.5mm}
\begin{enumerate}
\item [(1)]{\it $R$ is locally stable.}
\vspace{-.5mm}
\item [(2)]{\it $aR+bR=R$ with $a,b\in R\Longrightarrow~ \exists y\in R$ such that $R/(a+by)$ has stable range 1.}
\end{enumerate}
\vspace{-.5mm} {\it Proof.}\ \  This is an immediate consequence of Lemma 2.2.\hfill$\Box$

\vskip4mm Following McGovern, a ring $R$ has almost stable range 1 provided that every proper homomorphic image of $R$ has stable range 1 (cf. [4]). Every ring having stable range 1 has almost stale range 1, but the converse is not true, e.g., ${\Bbb Z}$.

\vskip4mm \hspace{-1.8em} {\bf Corollary 2.4.}\ \ {\it Every ring having almost stable range 1 is locally stable.}
\vskip2mm\hspace{-1.8em} {\it Proof.}\ \ Let $R$ have almost stable range 1. Suppose that $aR+bR=R$ with $a,b\in R$. If there exists a $y\in R$ such that $a+by\neq 0$. Then $R/(a+by)R$ has stable range 1. If $a+bz=0$ for all $z=0$, then $a=0$, and so $b\in U(R)$. Hence, $R/(a+b)R$ has stable range 1.
According to Theorem 2.3, $R$ is locally stable, as desired.\hfill$\Box$

\vskip4mm \hspace{-1.8em} {\bf Corollary 2.5.}\ \ {\it Every B$\acute{e}$zout ring having adequate range 1 is locally stable.}
\vskip2mm\hspace{-1.8em} {\it Proof.}\ \ \ Let $R$ be a B$\acute{e}$zout ring having adequate range 1. Suppose that $aR+bR=R$ for some $a,b\in R$. By hypothesis, there exists a $y\in R$ such that $a+by\in R$ is adequate. By virtue of [9, Theorem 2], $\overline{0}\in R/(a+by)R$ is adequate.
It follows from [9, Theorem 3] that $R/(a+by)R$ has stable range 1. Then by Theorem 2.3. $R$ is locally stable.\hfill$\Box$

\vskip4mm A
B$\acute{e}$zout ring in which every nonzero element is adequate
is called an adequate ring [7]. A ring $R$ is clean if every element in $R$ is the sum of an idempotent and
a unit. A ring $R$ is a neat ring if every proper homomorphic image is clean [3].

\vskip4mm \hspace{-1.8em} {\bf Example 2.6.}
\begin{enumerate}
\item [(1)] Every ring of stable range 1 is locally stable.
\item [(2)] Every adequate ring is locally stable.
\item [(3)] Every neat ring is locally stable.
\end{enumerate}
\vspace{-.5mm} {\it Proof.}\ \ (1) Every ring having stable range 1 has almost stable range 1. Then proving (1) by Corollary 2.4.

(2) As every adequate ring has adequate range 1, then by Corollary 2.5. it is locally stable.

(3) Let $aR+bR=R$ with $a,b\in R$. As in the proof of Corollary 2.4, there exists a $y\in R$ such that $a+by\neq 0$, so $R/(a+by)R$ is a proper homomorphic image of $R$; hence that it is clean. In view of [1, Corollary 1.3.15], $R/(a+by)$ has stable range 1. Therefore $R$ is locally stable, by Theorem 2.3.\hfill$\Box$

\vskip4mm Recall that a ring $R$ has stable range 2 provided that $aR+bR+cR=R$ with $a,b,c\in R\Longrightarrow ~\exists y,z\in R$ such that $(a+cy)R+(b+cz)R=R$.
\vskip4mm \hspace{-1.8em} {\bf Lemma 2.7.}\ \ {\it Every locally stable ring has stable range 2.}
\vskip2mm\hspace{-1.8em} {\it Proof.}\ \ Let $R$ be a locally stable ring. Suppose that $aR+bR+cR=R$ with $a,b,c\in R$.
Then there exist $y,z\in R$ such that $w:=a+by+cz$ and $R_w$ has stable range 1. Write $ak+bp+cq=1$ with $k,p,q\in R$. Then
$(a+by+cz)k+b(p-yk)+c(q-zk)=1$. Hence, $wR+bR+cR=1$, and so $wR_w+bR_w+cR_w=R_w$.
In virtue of Lemma 2.1, $w\in J(R_w)$, and so
$bR_a+cR_w=R_w$. Thus, we can find some $d\in R_w$ such that
$b+cd\in U(R_w)$. Hence, we have a $u\in S_w$ and a $v\in R$ such that $wR+(bu+cv)R=R$. Write $ws+(bu+cv)t=1$ and $up+wq=1$. We infer that
$ws+\big(bu+cv(up+wq)\big)t=1$. It follows that $w(s+cvqt)+(b+cvp)ut=1$.
Hence, $(a+by+cz)(s+cvqt)+(b+cvp)ut=1$, and so
$$\big(a+(b+cvp)y+c(z-vpy)\big)(s+cvqt)+(b+cvp)ut=1.$$ Thus,
$$\big(a+c(z-vpy)\big)(s+cvqt)+(b+cvp)\big(ut+y(s+cvqt)\big)=1.$$
Hence, $\big(a+c(z-vpy)\big)R+(b+cvp)R=R$. Therefore $R$ has stable range 2.
\hfill$\Box$

\vskip4mm \hspace{-1.8em} {\bf Lemma 2.8 [5, Theorem 2.5].}\ \ {\it Let $R$ be a Hermite ring. Then the following are
equivalent:}
 \vspace{-.5mm}
\begin{enumerate}
\item [(1)] {\it $R$ is an elementary divisor ring.} \vspace{-.5mm}
\item [(2)] {\it Every matrix $\left(
\begin{array}{cc}
a&0\\
b&c
\end{array}
\right)\in M_2(R)$ with $aR+bR+cR=R$ admits an elementary reduction.}
\end{enumerate}

\vskip4mm \hspace{-1.8em} {\bf Theorem 2.9.}\ \ {\it  Let $R$ be a locally stable ring. Then $R$ is an elementary divisor ring if and only if $R$ is a B$\acute{e}$zout ring.}
\vskip2mm\hspace{-1.8em} {\it Proof.}\ \  One direction is obvious.

Conversely, assume that $R$ is a B$\acute{e}$zout ring. In view of Lemma 2.7, $R$ has stable range 2. It follows from [7, Theorem 2.1.2] that $R$ is a Hermite ring. Suppose that
$A=\left(
\begin{array}{cc}
a&0\\
b&c
\end{array}
\right)\in M_2(R)$ with $aR+bR+cR=R$. We will suffice to prove that $A$ admits an elementary reduction.

Write $ax+by+cz=1$. Then $bR+(ax+cz)R=R$. As $R$ is locally stable, there exists some $t\in R$ such that $v=b+(ax+cz)t$ and $R_v$ has stable range 1. We have $$\left(
\begin{array}{cc}
1&0\\
xt&1
\end{array}
\right) \left(
\begin{array}{cc}
a&0\\
b&c
\end{array}
\right)\ \left(
\begin{array}{cc}
1&0\\
zt&1
\end{array}
\right)=\left(
\begin{array}{cc}
a&0\\
v&c
\end{array}
\right).$$ Since $R$ is a Hermite ring, there exists $P\in GL_2(R)$ such that $(v,c)P=(0,c^{\prime})$, and so $\left(
\begin{array}{cc}
a&0\\
v&c
\end{array}
\right)P=\left(
\begin{array}{cc}
a^{\prime}&b^{\prime}\\
0&c^{\prime}
\end{array}
\right)$. It is easily seen that $vR\subseteq c^{\prime}R$ and $a^{\prime}R+b^{\prime}R+c^{\prime}R=R$. Since $R_v$ has stable range 1, it follows by Lemma 2.2 that $R/vR$ has stable range 1, and so $R/c^{\prime}R$ has stable range 1, as
$R/c'R\cong (R/vR)/(c'R/vR)$. Also $\overline{a^{\prime}}(R/c^{\prime}R) + \overline{b^{\prime}}(R/c^{\prime}R)=R/c^{\prime}R$, hence there exists some $w\in R$ such that $\overline{b^{\prime}+a^{\prime}w}\in U(R/c^{\prime}R)$, and so $(\overline{b^{\prime}+a^{\prime}w)p}=\overline{1}$, for some $p\in R$, and then $(b'+aw)p+c'q=1$. One easily checks that
$$\begin{array}{ll}
&\left(
\begin{array}{cc}
&1\\
1&
\end{array}
\right)\left(
\begin{array}{cc}
c'&-(b'+aw)\\
p&q
\end{array}
\right)\left(
\begin{array}{cc}
a^{\prime}&b'\\
&c^{\prime}
\end{array}
\right)\left(
\begin{array}{cc}
1&w\\
&1
\end{array}
\right)\left(
\begin{array}{cc}
1&\\
-pa'&1
\end{array}
\right)\left(
\begin{array}{cc}
&1\\
1&
\end{array}
\right)\\
=&\left(
\begin{array}{cc}
1&\\
&a'c'
\end{array}
\right).
\end{array}$$
As $det\left(
\begin{array}{cc}
c'&-(b'+aw)\\
p&q
\end{array}
\right)=1$, we see that $\left(
\begin{array}{cc}
c'&-(b'+aw)\\
p&q
\end{array}
\right)\in GL_2(R)$. Thus, we prove that
$\left(
\begin{array}{cc}
a^{\prime}&b^{\prime}\\
&c^{\prime}
\end{array}
\right)$ admits a diagonal reduction. Therefore $\left(
\begin{array}{cc}
a&0\\
b&c
\end{array}
\right)$ admits a diagonal reduction, and so $R$ is an elementary divisor ring.\hfill$\Box$

\vskip4mm \hspace{-1.8em} {\bf Corollary 2.10 [4, Theorem 3.7].}\ \ {\it  Let $R$ have almost stable range 1. Then $R$ is an elementary divisor ring if and only if $R$ is a B$\acute{e}$zout ring.}
\vskip2mm\hspace{-1.8em} {\it Proof.}\ \ In view of Corollary 2.4, $R$ is a locally stable ring. This completes the proof, in terms of Theorem 2.9.\hfill$\Box$

\vskip4mm Immediately, we prove the known result: every B$\acute{e}$zout ring having stable range 1 is an elementary divisor ring.

\vskip4mm \hspace{-1.8em} {\bf Corollary 2.11.}\ \ {\it Let $R$ be a ring in which every element not in $J(R)$ is adequate.
Then $R$ is an elementary divisor ring if and only if $R$ is a B$\acute{e}$zout ring.}
\vskip2mm\hspace{-1.8em} {\it Proof.}\ \  One direction is obvious. Conversely assume that $R$ is a B$\acute{e}$zout ring in which  every element not in $J(R)$ is adequate.  We claim that $R$ is locally stable. Suppose that $aR+bR=R$ with $a,b\in R$. Case 1. $a\in J(R)$, then $b\not\in J(R)$, since $aR+bR=R$. Hence, $a+b\not\in J(R)$, so it is adequate and therefore $R/(a+b\cdot 1)R$  has stable range 1, as proof in Corollary 2.5. If $a\not\in J(R)$, then $a+b\cdot 0\not\in J(R)$. Hence  $a+b\cdot 0$ is adequate and so $R/(a+b\cdot 0)R$ has stable range 1. Accordingly, $R$ is locally stable, and then by Theorem 2.9, $R$ is an elementary divisor ring.\hfill$\Box$

\vskip4mm Kaplansky proved that for an
adequate ring being a Hermite ring was equivalent to being an
elementary divisor ring. This was extended to rings with
zero-divisors by Gillman and Henriksen (cf. [7, Theorem 1.2.13]). As an immediate consequence, we deduce that every adequate ring is an elementary divisor ring. In [2], Domsha and Vasiunyk considered domain of adequate range 1. We now extend their main result [2, Theorem 14] to B$\acute{e}$zout rings with zero divisors.

\vskip4mm \hspace{-1.8em} {\bf Corollary 2.12.}\ \ {\it  Let $R$ have adequate range 1. Then $R$ is an elementary divisor ring if and only if $R$ is a B$\acute{e}$zout ring.}
\vskip2mm\hspace{-1.8em} {\it Proof.}\ \  One direction is obvious. Conversely, assume that $R$ is a B$\acute{e}$zout ring. In view of Corollary 2.5, $R$ is locally stable. Therefore we complete the proof, by Theorem 2.9.\hfill$\Box$

\vskip4mm \hspace{-1.8em} {\bf Corollary 2.13.}\ \ {\it Let $R$ be a B$\acute{e}$zout domain. Then the following are equivalent:}
\vspace{-.5mm}
\begin{enumerate}
\item [(1)]{\it $R$ is locally stable.}
\vspace{-.5mm}
\item [(2)]{\it $aR+bR=R$ with $a,b\in R\Longrightarrow \exists y\in R$ such that $R/(a+by)$ is clean.}
\end{enumerate}
\vspace{-.5mm} {\it Proof.}\ \ $(1)\Leftarrow (2)$ In view of Theorem 2.9, $R$ is an elementary divisor ring. Thus, proving (2) by [8, Theorem 33].

$(2)\Rightarrow (1)$ Suppose that $aR+bR=R$ with $a,b\in R$. Then there exists a $y\in R$ such that $R/(a+by)$ is clean. In view of [1, Corollary 1.3.15],
$R/(a+by)$ has stable range 1. Therefore $R$ is locally stable, by Theorem 2.3.\hfill$\Box$

\vskip10mm\bc{\bf 3. Elementary Operations}\ec

\vskip4mm Using elementary operations are crucial in elementary matrix reduction. There are elementary divisor rings over which there exist matrices that have no elementary matrix reduction only by elementary operations. In current section, we shall show that matrices over locally stable rings could be reduced by elementary row operations and column ones. We now derive

\vskip4mm \hspace{-1.8em} {\bf Theorem 3.1.}\ \ {\it Let $R$ be a locally stable B$\acute{e}$zout ring. Then
for any matrix $A$ of size $n$ we can find such unimodular matrices $P\in GE_n(R)$ and $Q\in GL_n(R)$, that
$$PAQ=\left(
\begin{array}{cccc}
\varepsilon_1&&&\\
&\varepsilon_2&&\\
&&\ddots&\\
&&&\varepsilon_n
\end{array}
\right),$$ where $\varepsilon_i$ is a divisor of $\varepsilon_{i+1} (1\leq i\leq n-1)$.}
\vskip2mm\hspace{-1.8em} {\it Proof.}\ \ In view of Theorem 2.9, $R$ is an elementary divisor ring. Thus, for any $n\times n$ matrix $A, n>2$ there exist $P,Q\in GE_{n}(R)$ such that $PAQ$ is a diagonal matrix by [7, Theorem 4.3.3]. So we need only consider the case $n=2$. Let $A=\left(
\begin{array}{cc}
a&b\\
c&d
\end{array}
\right)$. Then we have a $U\in GL_2(R)$ such that
$AU=\left(
\begin{array}{cc}
x&0\\
y&z
\end{array}
\right).$
Clearly, $R$ is a B\'{e}zout ring. Write $xR+yR+zR=hR$. Then $x=x'h,y=y'h,z=z'h$ and $xp+yq+zk=h$. Set
$w=x'p+y'q+z'k-1$. Then $hw=0$ and $AU=h\left(
\begin{array}{cc}
x'&w\\
y'&z'
\end{array}
\right)$ with $x'R+y'R+z'R+wR=R$.
Similarly, we can find a $V\in GL_2(R)$ such that $(AU)V=h\left(
\begin{array}{cc}
a'&0\\
c'&d'
\end{array}
\right)$ with $a'R+c'R+d'R=R$. Write $a'x'+c'y'+d'z'=1$. Then $(-c')(-y)+(-a')(-x')+(-d')(-z')=1$
By hypothesis,
there exists some $p\in R$ such that $q:=-(c'+a'x'p+d'z'p)\in R$ and $R_q$ has stable range 1.
Thus,$$\left(
\begin{array}{cc}
&-1\\
1&
\end{array}
\right)\left(
\begin{array}{cc}
1&0\\
x'p&1
\end{array}
\right)
\left(
\begin{array}{cc}
a'&0\\
c'&d'
\end{array}
\right)\left(
\begin{array}{cc}
1&0\\
y'p&1
\end{array}
\right)=\left(
\begin{array}{cc}
q&-d'\\
a'&0
\end{array}
\right).$$
As $R$ is a Hermite ring, we have a $W\in GL_2(R)$ such that $(q,-d')W=(\alpha,0)$ for some $\alpha\in R$.

As $R_{q}$ has stable range 1. In virtue of Lemma 2.2, $R/qR$ has stable range 1.
Set $W^{-1}=(q_{ij})$. Then $q=\alpha q_{11}$. Clearly, $R/\alpha R\cong R/qR/\alpha R/qR$. It follows that $R/\alpha R$ has stable range 1.
Write $(a',0)W=(\gamma,\delta)$. Then
$\left(
\begin{array}{cc}
q&d'\\
a'&0
\end{array}
\right)W=\left(
\begin{array}{cc}
\alpha&\\
\gamma&\delta
\end{array}
\right)$, where $R/\alpha R$ has stable range 1 and $\alpha R+\gamma R+\delta R=R$.
As in the proof of Theorem 2.9,
we have a $P\in GE_2(R)$ and a $K\in GL_2(R)$ such that $$P\left(
\begin{array}{cc}
\alpha&\\
\gamma&\delta
\end{array}
\right)K=\left(
\begin{array}{cc}
1&\\
&\omega
\end{array}
\right).$$ Therefore $$P\left(
\begin{array}{cc}
&-1\\
1&
\end{array}
\right)\left(
\begin{array}{cc}
1&0\\
x'p&1
\end{array}
\right)
(AU)V\left(
\begin{array}{cc}
1&0\\
y'p&1
\end{array}
\right)WK=\left(
\begin{array}{cc}
h&\\
&h\omega
\end{array}
\right).$$ One easily checks that $$\left(
\begin{array}{cc}
&-1\\
1&
\end{array}
\right)=\left(
\begin{array}{cc}
1&-1\\
0&1
\end{array}
\right)\left(
\begin{array}{cc}
1&0\\
1&1
\end{array}
\right)\left(
\begin{array}{cc}
1&-1\\
0&1
\end{array}
\right)\in GE_2(R),$$ and therefore we complete the proof.
\hfill$\Box$

\vskip4mm \hspace{-1.8em} {\bf Corollary 3.2.}\ \ {\it Let $R$ be a B$\acute{e}$zout ring. If $R$ has almost stable range 1 or adequate range 1, then for any matrix of size $n$ we can find such unimodular matrices $P\in GE_n(R)$ and $Q\in GL_n(R)$, that
$$PAQ=\left(
\begin{array}{cccc}
\varepsilon_1&&&\\
&\varepsilon_2&&\\
&&\ddots&\\
&&&\varepsilon_n
\end{array}
\right),$$ where $\varepsilon_i$ is a divisor of $\varepsilon_{i+1} (1\leq i\leq n-1)$.}
\vskip2mm\hspace{-1.8em} {\it Proof.}\ \ In view of Corollary 2.4 and Corollary 2.5, $R$ is locally stable. Therefore the result follows, by Theorem 3.1.\hfill$\Box$

\vskip4mm \hspace{-1.8em} {\bf Theorem 3.3.}\ \ {\it  Let $R$ be a B$\acute{e}$zout ring. If $R_a$ has stable range 1 for all $a\not\in J(R)$, then
for any matrix $A$ of size $n$ we can find such unimodular matrices $P\in GE_n(R)$ and $Q\in GL_n(R)$, that
$$PAQ=\left(
\begin{array}{cccc}
\varepsilon_1&&&\\
&\varepsilon_2&&\\
&&\ddots&\\
&&&\varepsilon_n
\end{array}
\right),$$ where $\varepsilon_i$ is a divisor of $\varepsilon_{i+1} (1\leq i\leq n-1)$.}
\vskip2mm\hspace{-1.8em} {\it Proof.}\ \  Suppose that $aR+bR=R$ with $a,b\in R$. Case 1. $a\in J(R)$, then $b\not\in J(R)$. Hence, $a+b\not\in J(R)$. If $a\not\in J(R)$, then $a+b\cdot 0\not\in J(R)$. In any case, we can find a $y\in R$ such that $a+by\not\in J(R)$. By hypothesis, $R_{a+by}$ has stable range 1.
This implies that $R$ is locally stable. Accordingly, we obtain the result, by Theorem 3.1.\hfill$\Box$

\vskip4mm \hspace{-1.8em} {\bf Corollary 3.4.}\ \ {\it Let $R$ be a B$\acute{e}$zout ring. If $R_a$ is an adequate ring for all $a\not\in J(R)$, then for any matrix of size $n$ we can find such unimodular matrices $P\in GE_n(R)$ and $Q\in GL_n(R)$, that
$$PAQ=\left(
\begin{array}{cccc}
\varepsilon_1&&&\\
&\varepsilon_2&&\\
&&\ddots&\\
&&&\varepsilon_n
\end{array}
\right),$$ where $\varepsilon_i$ is a divisor of $\varepsilon_{i+1} (1\leq i\leq n-1)$.}
\vskip2mm\hspace{-1.8em} {\it Proof.}\ \ Let $a\not\in J(R)$. Given $bR_a+cR_a=R_a$ with $b,c\in R_a$, then $b(R_a/aR_a)+c(R_a/aR_a)=R_a/aR_a$. Since $R_a$ is adequate, it follows from $a\in R_a$ is adequate. In view of [9, Theorem 2], $\overline{0}\in R_a/aR_a$ is adequate. It follows by [9, Theorem 3] that
$R_a/aR_a$ has stable range 1. By virtue of Lemma 2.1, $a\in J(R_a)$, and so $R_a$ has stable range 1.
This completes the proof, by Theorem 3.3.\hfill$\Box$

\vskip4mm A ring $R$ is a VNL ring if for any $a\in R$, $a$ or $1-a$ is regular. Clearly, every regular ring $R$ (i.e., for any $a\in R$ there exists some $x\in R$ such that $a=axa$) is a VNL ring.

\vskip4mm \hspace{-1.8em} {\bf Corollary 3.5.}\ \ {\it Let $R$ be a B$\acute{e}$zout ring. If $R/aR$ is VNL for all $a\not\in J(R)$, then
for any matrix $A$ of size $n$ we can find such unimodular matrices $P\in GE_n(R)$ and $Q\in GL_n(R)$, that
$$PAQ=\left(
\begin{array}{cccc}
\varepsilon_1&&&\\
&\varepsilon_2&&\\
&&\ddots&\\
&&&\varepsilon_n
\end{array}
\right),$$ where $\varepsilon_i$ is a divisor of $\varepsilon_{i+1} (1\leq i\leq n-1)$.}
\vskip2mm\hspace{-1.8em} {\it Proof.}\ \ Let $a\not\in J(R)$. By hypothesis, $R/aR$ is a VNL ring. Let $x\in R/aR$. Then $x$ or $\overline{1}-x\in R/aR$ is regular. If $x$ is regular, we have a $y\in R/aR$ such that $x=xyx$. Set $z=yxy$. Then $x=xzx$ and $z=zxz$. Let $e=xz$ and $u=1-xz+x$. Then $e=e^2\in R,
u^{-1}=1-xz+z$ and $x=eu$. One easily checks that $\big(x-(1-e)\big)\big(u(1-e)-e\big)=u(2e-1)\in U(R)$. Hence, $x-(1-e)\in U(R)$.
This implies that $x\in R/aR$ is clean.
If $\overline{1}-x\in R/aR$, analogously, we can find an idempotent $f\in R/aR$ and a unit $v\in R/aR$ such that $\overline{1}-x=f+v$. Hence, $x=\big(\overline{1}-f\big)-v$ is clean. Thus, $R/aR$ is a clean ring. By virtue of [1, Corollary 1.3.15], $R/aR$ has stable range one. It follows by Lemma 2.2 that $R_a$ has stable range 1. Therefore we obtain the result, in terms of Theorem 3.3.\hfill$\Box$

\vskip4mm Let $R$ be a B$\acute{e}$zout VNL ring (e.g., regular ring). Then for any matrix of size $n$ we can find such unimodular matrices $P\in GE_n(R)$ and $Q\in GL_n(R)$, that
$$PAQ=\left(
\begin{array}{cccc}
\varepsilon_1&&&\\
&\varepsilon_2&&\\
&&\ddots&\\
&&&\varepsilon_n
\end{array}
\right),$$ where $\varepsilon_i$ is a divisor of $\varepsilon_{i+1} (1\leq i\leq n-1)$. This is an immediate consequence of Corollary 3.5.

\vskip4mm \hspace{-1.8em} {\bf Theorem 3.6.}\ \ {\it Let $R$ be a B$\acute{e}$zout ring. If $R_a$ or $R_{1-a}$ has stable range 1 for all $a\in R$, then for any matrix of size $n$ we can find such unimodular matrices $P\in GE_n(R)$ and $Q\in GL_n(R)$, that
$$PAQ=\left(
\begin{array}{cccc}
\varepsilon_1&&&\\
&\varepsilon_2&&\\
&&\ddots&\\
&&&\varepsilon_n
\end{array}
\right),$$ where $\varepsilon_i$ is a divisor of $\varepsilon_{i+1} (1\leq i\leq n-1)$.}
\vskip2mm\hspace{-1.8em} {\it Proof.}\ \ Given $aR+bR=R$ with $a,b\in R$, then $ax+by=1$ for some $x,y\in R$. Hence, $a(x-y)+(a+b)y=1$. By hypothesis, $R_{a(x+y)}$ or $R_{(a+b)y}$ has stable range 1.

Case 1. $R_{a(x+y)}$  has stable range 1.  In view of Lemma 2.2, $R/(a(x+y))R$ has stable range 1. Clearly, $R/aR\cong R/(a(x+y))R/aR/(a(x+y))R$. It follows that $R_a/aR$ has stable range 1. By using Lemma 2.2 again, $R_a$ has stable range 1. That is, $R_{a+b\cdot 0}$ has stable range 1.

Case 2. $R_{(a+b)y}$ has stable range 1. As in the proof in Case I, we see that $R_{a+b\cdot 1}$ has stable range 1.

Accordingly, we have a $z\in R$ such hat $R_{a+bz}$ has stable range 1. Therefore, $R$ is locally stable. This completes the proof, by
Theorem 3.1.\hfill$\Box$

\vskip4mm Following Domsha and Vasiunyk [2], a ring $R$ is a local adequate ring if for every $a\in R$, $a$ or $1-a$ is adequate. We are now ready to prove the following.

\vskip4mm \hspace{-1.8em} {\bf Corollary 3.7.}\ \ {\it Let $R$ be a B$\acute{e}$zout ring. If $R$ is a local adequate ring, then for any matrix of size $n$ we can find such unimodular matrices $P\in GE_n(R)$ and $Q\in GL_n(R)$, that
$$PAQ=\left(
\begin{array}{cccc}
\varepsilon_1&&&\\
&\varepsilon_2&&\\
&&\ddots&\\
&&&\varepsilon_n
\end{array}
\right),$$ where $\varepsilon_i$ is a divisor of $\varepsilon_{i+1} (1\leq i\leq n-1)$.}
\vskip2mm\hspace{-1.8em} {\it Proof.}\ \ Let $a\in R$. Then $a$ or $1-a$ is adequate. In view of [9, Theorem 2], $\overline{0}\in R/aR$ or $\overline{0}\in R/(1-a)R$ is adequate. It follows by [9, Theorem 3] that
 $R/aR$ or $R(1-a)R$ has stable range 1. In terms of Lemma 2.2, we see that $R_a$ or $R_{1-a}$ has stable range 1. Therefore the proof is true, by Theorem 3.6.\hfill$\Box$

 \vskip4mm As an immediate consequence, we have

 \vskip4mm \hspace{-1.8em} {\bf Corollary 3.8 [7, Theorem 4.3.5].}\ \ Let $R$ be an adequate ring. Then for any matrix of size $n$ we can find such unimodular matrices $P\in GE_n(R)$ and $Q\in GL_n(R)$, that
$$PAQ=\left(
\begin{array}{cccc}
\varepsilon_1&&&\\
&\varepsilon_2&&\\
&&\ddots&\\
&&&\varepsilon_n
\end{array}
\right),$$ where $\varepsilon_i$ is a divisor of $\varepsilon_{i+1} (1\leq i\leq n-1)$.

\vskip10mm\bc{\bf 4. Some Localizations}\ec

\vskip4mm Set $T=\{ a~\mid~R_a~$ has stable range one $\}$. A subset $W$ of a ring $R$ is called saturated if it is closed under taking divisors: i.e.,
whenever a product $xy$ is in $S$, the elements $x$ and $y$ are in $W$ too. We now derive

\vskip4mm \hspace{-1.8em} {\bf Lemma 4.1.}\ \ {\it Let $R$ be a ring. Then $T$ be a saturated multiplicatively closed set of $R$.}
\vskip2mm\hspace{-1.8em} {\it Proof.}\ \ Clearly, $1\in T$. Let $a,b\in T$. By hypothesis, $R_a$ and $R_b$ has stable range 1. In view of Lemma 2.2, $R/aR$ and $R/bR$ has stable range 1. According to [1, Theorem 1.1.13], $R/(aR\cap bR)$ has stable range 1. Clearly, $$R/(aR\cap bR)\cong R/(ab)R/(aR\cap bR)/(ab)R;$$ and so $R/(ab)R/(aR\cap bR)/(ab)R$ has stable range 1. Since $\big((aR\cap bR)/(ab)R\big)^2=0$, we see that
$(aR\cap bR)/(ab)R\subseteq J\big(R/(ab)R$. This implies that $R/(ab)R$ has stable range 1. That is, $R_{ab}$ has stable range 1, by Lemma 2.2. Hence, $ab\in T$.

Assume that $c=de$ with $c\in T$. Then $R/cR$ has stable range 1. As $cR\subseteq dR$, we see that $R/dR\cong R/cR/dR/cR$, whence,
$R/dR$ has stable range 1. Thus, $d\in T$, as desired.\hfill$\Box$

\vskip4mm \hspace{-1.8em} {\bf Lemma 4.2.}\ \ \ \ {\it If $R$ is a B$\acute{e}$zout ring, then so is $R_T$.}
\vskip2mm\hspace{-1.8em} {\it Proof.}\ \ Let $R$ be a B$\acute{e}$zout ring, and let $I=\frac{x}{t}T+\frac{y}{t}T$ where $x,y \in R, t\in T$. Since $R$ is  a B$\acute{e}$zout ring, there exists a $z\in R$ such that $xR+yR=zR$. Let $\frac{x}{t}\frac{r_{1}}{s_{1}}+\frac{y}{t}\frac{r_{2}}{s_{2}}\in I$. Then there exists some $a\in R$ such that $xr_{1}s_{2}+yr_{2}s_{1}=az$, so we have $\frac{x}{t}\frac{r_{1}}{s_{1}}+\frac{y}{t}\frac{r_{2}}{s_{2}}=\frac{a}{s_1s_2}\frac{z}{t}\in \frac{z}{t}T$. Write $z=xa+yb$. For any $\frac{z}{t}\frac{r}{s}\in \frac{z}{t}T$, we have $\frac{z}{t}\frac{r}{s}=\frac{ar}{s}\frac{x}{t}+\frac{br}{s}\frac{y}{t}\in I$. Thus $I=\frac{z}{t}T$, as desired. \hfill$\Box$

\vskip4mm \hspace{-1.8em} {\bf Theorem 4.3.}\ \ {\it Let $R$ be a ring. If $R_T$ has adequate range 1, then $R$ is an elementary divisor ring if and only if $R$ is a Hermite ring.}
\vskip2mm\hspace{-1.8em} {\it Proof.}\ \ $\Longrightarrow$ This is obvious.

$\Longleftarrow$ Let $R$ be a Hermite ring. Then it is a B$\acute{e}$zout ring. In view of Lemma 4.2, $R_T$ is a B$\acute{e}$zout ring. In view of Corollary 2.10, $R_T$ is an elementary divisor ring. In light of Lemma 2.8, we will suffice to prove that every matrix $A=\left(
\begin{array}{cc}
a&0\\
b&c
\end{array}
\right)$ with $aR+bR+cR=R$ admits a diagonal reduction. As $aR_T+bR_T+cR_T=R_T$, it follows by Theorem 1.1 that we can find some $p,q\in R$ and $s_1,s_2\in T$ such that $$(aps_1^{-1}+bqs_2^{-1})R_T+cqs_2^{-1}R_T=R_T.$$
Thus, we have some $k,l,r,t\in R$ and $s\in T$ such that $$(ak+bl)r+clt=s.$$
Since $R$ is a Hermite ring, we can find some $k',l'\in R$ such that $k=k'x,l=l'x$ and $k'R+l'R=R$.
Hence, $x\big((ak'+bl')r+cl't\big)=s$. In view of Lemma 4.1, $(ak'+bl')r+cl't\in T$. Thus, we may assume that $kR+lR=R$.
Similarly, we can find some $r',t'\in R$ such that $r=r'y,t=t'y$ and $r'R+t'R=R$. Thus, $y\big((ak+bl)r'+clt'\big)=s$.
It follows that $(ak+bl)r'+clt'\in T$, by Lemma 4.1. Hence, we may assume that $rR+tR=R$.
Write $k\alpha+l\beta=1$ and $r\gamma+t\delta=1$. Then we check that
$$\left(
\begin{array}{cc}
&1\\
1&
\end{array}
\right)\left(
\begin{array}{cc}
\beta&-\alpha\\
k&l
\end{array}
\right)A\left(
\begin{array}{cc}
\delta&r\\
-\gamma&t
\end{array}
\right)\left(
\begin{array}{cc}
&1\\
1&
\end{array}
\right)=\left(
\begin{array}{cc}
s&*\\
**&*
\end{array}
\right).$$
Clearly, $\left(
\begin{array}{cc}
\beta&-\alpha\\
k&l
\end{array}
\right), \left(
\begin{array}{cc}
\delta&r\\
-\gamma&t
\end{array}
\right)\in GL_2(R)$, as their determinants are both $1$. Since $R$ is a Hermite ring, there exists a $Q\in GL_2(R)$ such that $(s,*)Q=(d,0)$. Write $Q^{-1}=(q_{ij})$. Then $s=dq_{11}$. In view of Lemma 4.1, $d\in T$. Therefore
$\left(
\begin{array}{cc}
s&*\\
**&*
\end{array}
\right)Q=\left(
\begin{array}{cc}
d&0\\
e&f
\end{array}
\right)$ for some $e,f\in R$. One easily checks that $dR+eR+fR=R$. In light of Lemma 2.2, $R/dR$ has stable range 1.
Clearly, $\overline{e}(R/dR)+\overline{f}(R/dR)=R/dR$. Then, we can find a $w\in R$ such that $\overline{e+fw}\in U(R/dR)$. It follows that
$dR+(e+fw)R=R$. Write $dm+(e+fw)n=1$ for some $m,n\in R$. Then we verify that
$$\left(
\begin{array}{cc}
m&n\\
-(e+fw)&d
\end{array}
\right)\left(
\begin{array}{cc}
d&0\\
e&f
\end{array}
\right)\left(
\begin{array}{cc}
1&\\
w&1
\end{array}
\right)\left(
\begin{array}{cc}
1&-nf\\
&1
\end{array}
\right)=\left(
\begin{array}{cc}
1&\\
&df
\end{array}
\right).$$ Obviously, $\left(
\begin{array}{cc}
m&n\\
-(e+fw)&d
\end{array}
\right)\in GL_2(R)$. Therefore $A$ admits an elementary diagonal reduction. This completes the proof.\hfill$\Box$

\vskip4mm \hspace{-1.8em} {\bf Corollary 4.4.}\ \ {\it Let $R$ be a domain. If $R_T$ is an adequate ring, then $R$ is an elementary divisor ring if and only if $R$ is a B$\acute{e}$zout ring.}
\vskip2mm\hspace{-1.8em} {\it Proof.}\ \ One direction is obvious. Conversely, assume that $R$ is a B$\acute{e}$zout ring. Since $R$ is a domain, it is a Hermite ring. As every adequate ring has adequate range 1, we complete the proof, by Theorem 4.3.\hfill$\Box$

\vskip4mm Set $S=\{ a\in R~\mid~a\in R$ is adequate $\}$. As in the proof of [7, Proposition 3.2.2], we see that the product of two adequate elements is adequate. Thus, $S$ is a multiplicatively closed set of $R$. Zabavsky has ever studied matrix reduction on $R_S$ for a B$\acute{e}$zout domain $R$ [7, Theorem 3.2.2 ]. We now explore $R_S$ for a ring which may contain many zero divisors.

\vskip4mm \hspace{-1.8em} {\bf Theorem 4.5.}\ \ {\it Let $R$ be a ring. If $R_S$ has adequate range 1, then $R$ is an elementary divisor ring if and only if $R$ is a Hermite ring.}
\vskip2mm\hspace{-1.8em} {\it Proof.}\ \ $\Longrightarrow$ This is obvious.

$\Longleftarrow$ Let $R$ be a Hermite ring, for being an elementary divisor ring, in light of Lemma 2.8, we need only to prove that every matrix $ A=\left(
\begin{array}{cc}
a&0\\
b&c
\end{array}
\right)$ with $aR+bR+cR=R$ admits a diagonal reduction. We have $aR_S+bR_S+cR_S=R_S$. As $R$ is a B$\acute{e}$zout ring,  then so is $R_S$.
By virtue of Corollary 2.5, $R_S$ is locally stable. Thus, by Theorem 2.9, $R_S$ is an elementary divisor ring. Now by using Theorem 1.1, we can find some elements $p,q\in R$ and $s_1,s_2\in S$ such that $(aps_1^{-1}+bqs_2^{-1})R_S+cqs_2^{-1}R_S=R_S.$ As in the proof of Theorem 4.3, there is an equivalent matrix $ B=\left(
\begin{array}{cc}
z&0\\
x&y
\end{array}
\right)$ for $A$ with $xR+yR+zR=R$ and $z\in S$. Since $z\in R$ is an adequate element, then $R/zR$ has stable range 1.
 Also $\overline{x}(R/zR) + \overline{y}(R/zR)=R/zR$, hence there exists some $v\in R$ such that $\overline{x+yv}\in U(R/zR)$. So $(\overline{x+yv)u}=\overline{1}$, for some $u\in R$, and then $(x+yv)u+zt=1$ for some $t\in R$. It is obvious that
$$
\left(
\begin{array}{cc}
t&u\\
-(x+yv)&z
\end{array}
\right)\left(
\begin{array}{cc}
z&\\
x&y
\end{array}
\right)\left(
\begin{array}{cc}
1&\\
v&1
\end{array}
\right)\left(
\begin{array}{cc}
1&-uy\\
&1
\end{array}
\right)=\left(
\begin{array}{cc}
1&\\
&yz
\end{array}
\right).$$
As $det\left(
\begin{array}{cc}
t&u\\
-(x+yv)&z
\end{array}
\right)=1$, we see that $\left(
\begin{array}{cc}
t&u\\
-(x+yv)&z
\end{array}
\right)\in GL_2(R)$. Thus,
$\left(
\begin{array}{cc}
z&\\
x&y
\end{array}
\right)$ admits an elementary diagonal reduction, and therefore $A$ admits an elementary diagonal reduction. So the theorem is true.\hfill$\Box$

\vskip4mm \hspace{-1.8em} {\bf Corollary 4.6.}\ \ {\it Let $R$ be a ring. If $R_S$ has stable range 1, then $R$ is an elementary divisor ring if and only if $R$ is a Hermite ring.}
\vskip2mm\hspace{-1.8em} {\it Proof.}\ \ Since $R_S$ has stable range 1, it has adequate range 1. Therefore we obtain the result, by Theorem 4.5.\hfill$\Box$

\vskip4mm \hspace{-1.8em} {\bf Corollary 4.7.}\ \ {\it Let $R$ be a ring. If $R_S$ is an adequate ring, then $R$ is an elementary divisor ring if and only if $R$ is a Hermite ring.}
\vskip2mm\hspace{-1.8em} {\it Proof.}\ \ As every adequate ring has adequate range 1, the corollary is established, in terms of
Theorem 4.5.\hfill$\Box$

\vskip15mm \bc{\bf REFERENCES}\ec \vskip4mm {\small
\re{1} H.
Chen, {\it Rings Related Stable Range Conditions}, Series in
Algebra 11, World Scientific, Hackensack, NJ, 2011.

\re{2} O.V. Domsha and I.S. Vasiunyk, Combining local and adequate rings, {\it Book of abstracts of the International Algebraic Conference}, Taras Shevchenko National University of Kyiv, Kyiv, Ukraine, 2014.

\re{3} W.W. McGovern, Neat rings, {\it J. Pure Appl. Algebra}, {\bf 205}(2006), 243--265.

\re{4} W.W. McGovern, B\'{e}zout rings with almost stable range
$1$, {\it J. Pure Appl. Algebra}, {\bf 212}(2008), 340--348.

\re{5} M. Roitman, The Kaplansky condition and rings of almost
stable range $1$, {\it Proc. Amer. Math. Soc.}, {\bf 141}(2013),
3013--3018.

\re{6} B.V. Zabavsky, Elementary reduction of matrices over adequate domain,
{\it Math. Stud.}, {\bf 17}(2002), 115--116.

\re{7} B.V. Zabavsky, Diagonal Reduction of Matrices over Rings,
Mathematical Studies Monograph Series, Vol. XVI, VNTL Publisher, 2012.

\re{8} B.V. Zabavsky, Diagonal reduction of matrices over finite stable range rings,
{\it Math. Stud.}, {\bf 41}(2014), 101-108.

\re{9} B.V. Zabavsky and S.I. Bilavska,
Every zero adequate ring is an exchange ring, {\it J. Math. Sci.},
{\bf 187}(2012), 153--156.

\re{10} B.V. Zabavsky and O. Domsha, Diagonalizability theorem for matrices
over certain domains, {\it Algebra $\&$ Discrete Math.}, {\bf 12}(2011), 132--139.

\end{document}